\documentclass[12pt,a4paper]{article}

\usepackage[T1]{fontenc}
\usepackage[latin2]{inputenc}
\usepackage{amsmath}
\usepackage{amssymb}
\usepackage{graphicx}
\newtheorem{thm}{Theorem}
\newtheorem{lem}{Lemma}

\newtheorem{cor}{Corollary}
\newtheorem{prob}{Problem}
\newtheorem{con}{Conjecture}
\begin{document}

\title{\vspace*{-3cm}On a problem of Chen and Fang related to infinite additive complements}
\author{S\'andor Z. Kiss \thanks{Institute of Mathematics, Budapest
University of Technology and Economics, H-1529 B.O. Box, Hungary;
kisspest@cs.elte.hu;
This author was supported by the National Research, Development and Innovation \
Office NKFIH Grant No. K115288 and K129335.
This paper was supported by the J\'anos Bolyai Research Scholarship of the Hungarian Academy of Sciences. Supported by the \'UNKP-20-5 New National Excellence Program
of the Ministry for Innovation and Technology from the source of the National Research, Development and Innovation Fund.}, Csaba S\'andor \thanks{Institute of Mathematics, Budapest University of
Technology and Economics, MTA-BME Lend\"ulet Arithmetic Combinatorics Research \
Group H-1529 B.O. Box, Hungary, csandor@math.bme.hu.
This author was supported by the NKFIH Grants No. K129335. Research supported by
the Lend\"ulet program of the Hungarian Academy of Sciences (MTA), under grant number LP2019-15/2019.}}

\date{}

\maketitle

\begin{abstract}
\noindent Two infinite sets $A$ and $B$ of nonnegative integers are called additive complements if their sumset contains every nonnegative integer. In 1964, Danzer constructed infinite additive complements $A$ and $B$ with $A(x)B(x) = (1 + o(1))x$ as $x \rightarrow \infty$, where $A(x)$ and $B(x)$ denote the counting function of the sets $A$ and $B$, respectively. In this paper we solve a problem of Chen and Fang by extending the construction of Danzer.


\end{abstract}

\section{Introduction}
Let $\mathbb{N}$ be the set of nonnegative integers and let $A$ and $B$ be infinite sets of nonnegative integers. We define their sum by $A + B = \{a+b: a\in A, b\in B\}$. We say $A$ and $B$ are infinite additive complements if their sum contains all nonnegative integers i.e., $A + B = \mathbb{N}$. Let $A(x)$ be the number of elements of $A$ up to $x$ i.e.,
\[
A(x) = \sum_{\overset{a\in A}{a \le x}}1.
\]
Since $A$ and $B$ are infinite additive complements, every nonnegative integer $x$ can be written in the form $a + b = x$, where
$a\in A$, $b \in B$. Then clearly \cite{sarkoszeme} we have $A(x)B(x) \ge x + 1$, which implies that
\[
\limsup_{x \rightarrow \infty}\frac{A(x)B(x)}{x} \ge \liminf_{x \rightarrow \infty}\frac{A(x)B(x)}{x} \ge 1.
\]
According to a conjecture of H. Hanani \cite{erd}, the above result can be sharpened in the following way.
\begin{con}[Hanani, 1957]
If $A$ and $B$ are infinite additive complements, then
\[
\limsup_{x \rightarrow \infty}\frac{A(x)B(x)}{x} > 1.
\]
\end{con}
Later, Danzer \cite{Danzer} disproved the above conjecture of Hanani.

\begin{thm}[Danzer, 1964]
There exist infinite additive complements $A$ and $B$ such that
\[
\lim_{x \rightarrow \infty}\frac{A(x)B(x)}{x} = 1.
\]
\end{thm}

Let $A_{1}, \dots{} ,A_{r}$ be infinite sets of nonnegative integers. We define their sum by
$A_{1} + A_{2} + \dots{} + A_{r} = \{a_{1} + a_{2} + \dots{} + a_{r}: a_{i}\in A_{i}, 1 \le i \le r\}$.
Chen and Fang extended the notion of additive complements to more than two sets in the following way \cite{Chen}. The infinite sets $A_{1}, \dots{} ,A_{r}$ of nonnegative integers are said to form infinite additive complements if their sum contains all nonnegative integers. Again, it is easy to see that
$A_{1}(x)\cdots{} A_{r}(x) \ge (A_{1} + \dots{} + A_{r})(x) = x + 1$, thus
\[
\liminf_{x \rightarrow \infty}\frac{A_{1}(x)\cdots{} A_{r}(x)}{x} \ge 1.
\]
Furthermore, they posed the following problem.
\begin{prob}
For each integer $r \ge 3$ find additive complements $A_{1}, \dots{} ,A_{r}$ such that
\[
\lim_{x \rightarrow \infty}\frac{A_{1}(x)\cdots{} A_{r}(x)}{x} = 1.
\]
\end{prob}
In this paper we solve this problem. Note that our construction is the extension of Danzer's result to $r > 2$.

\begin{thm}
For each integer $h \ge 2$ there exist infinite sets of nonnegative integers
$A_{1}, \dots{} ,A_{h}$ with the following properties:
\begin{itemize}
    \item[(1)] $A_{1} + \dots{} + A_{h} = \mathbb{N}$,
    \item[(2)] $A_{1}(x) \cdots{}  A_{h}(x) = (1 + o(1))x$ as $x \rightarrow \infty$.
\end{itemize}
\end{thm}
Let $R_{A+B}(n)$ be the number of representations of the integer $n$ in the form $a + b = n$, where $a\in A$, $b\in B$. W. Narkiewicz \cite{Nar} proved the following theorem.

\begin{thm}[Narkiewicz, 1960]
If $R_{A+B}(n) \ge C$ for every sufficiently large integer $n$, where $C$ is a constant and
\[
\limsup_{x \rightarrow \infty}\frac{A(x)B(x)}{x} \le C,
\]
then
\[
\lim_{x \rightarrow \infty}\frac{A(2x)}{x} = 1,
\]
or
\[
\lim_{x \rightarrow \infty}\frac{B(2x)}{x} = 1.
\]
\end{thm}
Additive complements $A$, $B$ are called exact if $A(x)B(x) = (1 + o(1))x$ as $x \rightarrow \infty$.
For any $h \ge 2$ integer let us define the system of sets $\mathcal{A}_{h}$ by
\[
\mathcal{A}_{h}= \{A \subset \mathbb{N}: \text{there exist} \hspace*{1mm} A_{2}, \dots{} ,A_{h} \subset \mathbb{N},
\]
\[
A + A_{2} + \dots{} + A_{h} = \mathbb{N}, A(x)\cdot A_{2}(x) \cdots{}  A_{h}(x) = (1 + o(1))x \hspace*{1mm} \text{as} \hspace*{1mm} x \rightarrow \infty\}.
\]
Theorem 3 implies that $\mathcal{A}_{h} \ne \emptyset$ for every  $h \ge 2$. We prove that the $\mathcal{A}_{h}$'s form an infinite chain.

\begin{thm}
We have $\mathcal{A}_{2} \supseteq \mathcal{A}_{3} \supseteq \dots{}$
\end{thm}
According to Theorem 3,
if $A \in \mathcal{A}_{2}$, then $A(x) = x^{o(1)}$ or $A(x) = x^{1+o(1)}$ as $x \rightarrow \infty$. Then for any $h \ge 2$, $A \in \mathcal{A}_{h}$ implies that $A(x) = x^{o(1)}$ or $A(x) = x^{1+o(1)}$ as
$x \rightarrow \infty$. If the sets $A_{1}, \dots{} ,A_{h} \subset \mathbb{N}$ satisfy $A_{1} + \dots{} + A_{h} = \mathbb{N}$ and $A_{1}(x) \cdots{}  A_{h}(x) = (1 + o(1))x$ as $x \rightarrow \infty$, then $A_{i}(x) = x^{1+o(1)}$ or $A_{i}(x) = x^{o(1)}$ for every $1 \le i \le h$ while $x \rightarrow \infty$.
As a corollary, one can get from Theorem 4 that

\begin{cor}
Let $A_{1}, \dots{} ,A_{h}$ be infinite sets of nonnegative integers such that $A_{1} + \dots{} + A_{h} = \mathbb{N}$ and
\[
A_{1}(x) \cdots{}  A_{h}(x) = (1 + o(1))x
\]
as $x \rightarrow \infty$. Then there exists an index $i$ such that $A_{i}(x) = x^{1+o(1)}$ and $A_{j}(x) = x^{o(1)}$ for every $1 \le j \le h$ with $j \ne i$ as $x \rightarrow \infty$.
\end{cor}
We pose the following problems for further research.

\begin{prob}
Does $\mathcal{A}_{h} \ne \mathcal{A}_{h+1}$ hold for every $h\ge 2$?
\end{prob}

\begin{prob}
Assume that $A_{1} + \dots{} + A_{h} = \mathbb{N}$ and $A_{1}(x) \cdots{}  A_{h}(x) = (1 + o(1))x$ hold as $x \rightarrow \infty$.
Does there exist a permutation $i_{1}, \dots{} ,i_{h}$ of the indices $1, \dots{} ,h$ such that $A_{i_{j}}(x) = (A_{i_{j-1}}(x))^{o(1)}$ for every $2 \le j \le h$ as $x \rightarrow \infty$?
\end{prob}
The statement in Problem 2 holds for $h = 2$.

The exact complemets have been investigated by many authors in the last few decades. In particular, they studied what kind of sets $A$ of nonnegative integers with $A(x) = x^{o(1)}$ as $x \rightarrow \infty$ have exact additive complement. It was proved in \cite{Danzer} that the sequence $a_{n} = (n!)^{2} + 1$ has an exact complement. In \cite{Ruzsa2} Ruzsa showed that the set of the powers of an integer $a \ge 3$ has an exact complement.
Furthermore, in \cite{Ruzsa1} he proved that the set of powers of $2$ has an exact complement. Moreover, he also proved in \cite{Ruzsa1} that $A = \{a_{1}, a_{2},\dots{}\}$ with $1 \le a_{1} < a_{2} < \dots{}$ has an exact complement if $\lim_{n \rightarrow \infty}\frac{a_{n+1}}{na_{n}} = \infty$. In view of these results, it is natural to ask
\begin{prob}
Is it true that if $A \in \mathcal{A}_{2}$, $A(x) = x^{o(1)}$ as $x \rightarrow \infty$, then $A(x) = O(\log x)$?
\end{prob}

\section{Proof of Theorem 2}
For any nonnegative integers $a < b$, let us define $[a,b] = \{x\in \mathbb{N}: a\le x \le b\}$.
For any nonnegative integer $k$, we define $[a,b]\cdot k = \{kx: x\in [a,b] \cap \mathbb{N}\}$.
The following lemma plays the key role in the proof of Theorem 2.
\begin{lem}
Assume that
$A_{1}, \dots{} ,A_{h} \subset \mathbb{N}$ infinite subsets with the following properties:
\begin{itemize}
    \item[(1)] $A_{1} + \dots{} + A_{h} = \mathbb{N}$,
    \item[(2)] there exists a monotone increasing arithmetic function $f_{h}(n) \ge 0$ with $$\lim_{n\rightarrow \infty}f_{h}(n) = \infty$$ such that the equation
    $a_{1} + \dots{} + a_{h} = n$, $a_{i} \in A_{i}$ has a solution with
    $a_{i} \ge f_{h}(n)$,
    \item[(3)] $A_{1}(x) \cdots{}  A_{h}(x) = (1 + o(1))x$ as $x \rightarrow \infty$.
\end{itemize}
For $m \in \mathbb{N}$, let $g(m)$ be a strictly increasing function such that $g(f_{h}(n)) \ge n^{2}$ for every $n \in \mathbb{N}$. Furthermore, for $1 \le i \le h$, let $B_{i} = \{g(a)! + a: a\in A_{i}\}$ and define the sets of integers
\[
B_{h+1} = [0,g(6)!+h(g(6)-1)!-1]\cup
\]
\[
\left\{\bigcup_{n \ge 6}\left[\frac{g(n)!}{n-\lceil\sqrt{n}\rceil}-2,\frac{g(n+1)!+h(g(n+1)-1)!}{n-\lceil\sqrt{n}\rceil}\right]\cdot (n-\lceil\sqrt{n}\rceil)\right\}.
\]
Then
\begin{itemize}
    \item[(i)] $B_{1} + \dots{} + B_{h+1} = \mathbb{N}$,
    \item[(ii)] there exists a monotone increasing arithmetic function $f_{h+1}(n) \ge 0$ with $$\lim_{n\rightarrow \infty}f_{h+1}(n) = \infty$$ such that the equation
    $b_{1} + \dots{} + b_{h+1} = n$, $b_{i} \in B_{i}$ has a solution with
    $b_{i} \ge f_{h+1}(n)$,
    \item[(iii)] $B_{1}(x) \cdots{}  B_{h+1}(x) = (1 + o(1))x$ as $x \rightarrow \infty$.
\end{itemize}
\end{lem}

\subsection{Proof of the lemma}

Now we prove that for any $N \ge 6$,
\[
B_{1} + \dots{} + B_{h} + \left[\frac{g(N)!}{N-\lceil\sqrt{N}\rceil}-2,\frac{g(N+1)!+h(g(N+1)-1)!}{N-\lceil\sqrt{N}\rceil}\right]\cdot (N-\lceil\sqrt{N}\rceil)
\]
\[
\supseteq [g(N)! -2(N-\lceil\sqrt{N}\rceil) + h(g(N)-1)!+N, g(N+1)! +h(g(N+1)-1)!].
\]
Consider an element from the interval on the right hand side i.e., let $y$ be
\[
g(N)! -2(N-\lceil\sqrt{N}\rceil) + h(g(N)-1)!+N \le y \le g(N+1)! +h(g(N+1)-1)!.
\]
It is clear that there exists an $\lceil \sqrt{N} \rceil \le m \le N - 1$
with $y \equiv m \pmod{N-\lceil\sqrt{N}\rceil}$. By (2), there exist $a_{1}, \dots{} ,a_{h}$ integers with
$a_{i} \in A_{i}$ such that $m = a_{1} + \dots{} + a_{h}$ and $a_{i} \ge f_{h}(m)$. Since $f_{h}(m)$ is a monotone increasing function and $g(m)$ is a strictly increasing function, then we have
\[
g(a_{i}) \ge g(f_{h}(m)) \ge g(f_{h}(\lceil\sqrt{N}\rceil)) \ge (\lceil\sqrt{N}\rceil)^{2} \ge N
\]
and so $g(a_{i})! \equiv 0 \pmod{N-\lceil\sqrt{N}\rceil}$. Let $b_{i} = g(a_{i})! + a_{i}$. Then
$b_{i} \in B_{i}$ for every $1 \le i \le h$. It follows that
\[
\sum_{i=1}^{h}b_{i} = \sum_{i=1}^{h}(g(a_{i})! + a_{i}) \equiv \sum_{i=1}^{h}a_{i} \equiv m \equiv y \pmod{N-\lceil\sqrt{N}\rceil},
\]
which implies that $\frac{y-(b_{1} + \dots{} + b_{h})}{N-\lceil\sqrt{N}\rceil}$ is an integer and clearly
\[
y = b_{1} + \dots{} + b_{h} + \frac{y-(b_{1} + \dots{} + b_{h})}{N-\lceil\sqrt{N}\rceil}\cdot (N-\lceil\sqrt{N}\rceil).
\]
In view of these facts, it is enough to show that
\[
\frac{g(N)!}{N-\lceil\sqrt{N}\rceil}-2 \le \frac{y-(b_{1} + \dots{} + b_{h})}{N-\lceil\sqrt{N}\rceil} \le \frac{g(N+1)!+h(g(N+1)-1)!}{N-\lceil\sqrt{N}\rceil}.
\]
Since $g(n)$ is a strictly increasing function, we have
\[
0 \le b_{i} = g(a_{i})! + a_{i} \le g(m) + m \le g(N - 1)! + N - 1 < (g(N)-1)! + N
\]
and so
\[
0 \le \sum_{i=1}^{h}b_{i} < h((g(N) - 1)! + N).
\]
It follows that
\[
\frac{y-(b_{1} + \dots{} + b_{h})}{N-\lceil\sqrt{N}\rceil} \ge \frac{y-h((g(N) - 1)! +    N)}{N-\lceil\sqrt{N}\rceil}
\]
\[
\ge \frac{g(N)! -2(N-\lceil\sqrt{N}\rceil) + h((g(N)-1)!+N)-h(g(N)-1)!+N)}{N-\lceil\sqrt{N}\rceil}
\]
\[
= \frac{g(N)!}{N-\lceil\sqrt{N}\rceil} - 2
\]
and
\[
\frac{y-(b_{1} + \dots{} + b_{h})}{N-\lceil\sqrt{N}\rceil} \le \frac{y}{N-\lceil\sqrt{N}\rceil} \le \frac{g(N+1)! +h(g(N+1)-1)!}{N-\lceil\sqrt{N}\rceil}.
\]
Thus for $N \ge 6$, we have
\begin{align*}
&B_{1} + \dots{} + B_{h+1} \supseteq \\
&[g(N)! -2(N-\lceil\sqrt{N}\rceil) + h(g(N)-1)!+N, g(N+1)! +h(g(N+1)-1)!] \supseteq \\
&[g(N)! + h(g(N)-1)!, g(N+1)! + h(g(N+1)-1)!].
\end{align*}
This implies that
\[
B_{1} + \dots{} + B_{h+1} \supseteq \bigcup_{N \ge 6}[g(N)! + h(g(N)-1)!, g(N+1)! + h(g(N+1)-1)!]
\]
\[
= [g(6)! + h(g(6)-1)!, +\infty).
\]

Moreover, for $1 \le i \le h$, $0\in B_{i}$ and $B_{h+1} \supseteq [0, g(6)! + h(g(6)-1)!-1]$, thus
$[0, g(6)! + h(g(6)-1)!-1] \subseteq B_{1} + \dots{} + B_{h+1}$ and so
$B_{1} + \dots{} + B_{h+1} = \mathbb{N}$, which proves (i).

If $g(N)! + h(g(N)-1)! \le n \le g(N+1) + h(g(N+1)-1)!$, then there exists a representation $n = b_{1} + \dots{} + b_{h+1}$, where $b_{i} = g(a_{i}) + a_{i} \ge a_{i} \ge f_{h}(\lceil\sqrt{N}\rceil)$ and $b_{h+1} \ge g(N)! - 2(N-\lceil\sqrt{N}\rceil) \ge N! - 2(N-\lceil\sqrt{N}\rceil)$, which proves (ii) with a suitable function $f_{h+1}(n)$.

 To prove (iii) we assume that $g(N)! + h(g(N) - 1)! \le x \le g(N+1)! + h(g(N+1) - 1)!$. Since $g(N)$ is strictly increasing, $g(N+2h) \ge g(N+1) + h$. This implies that
 \[
 x \le g(N+1)! + h(g(N+1) - 1)! = (g(N+1)+h)(g(N+1)-1)!
 \]
 \[
 \le g(N+2h)(g(N+1)-1)! < g(N+2h)! + N + 2h
 \]
 and
 \[
 x \ge g(N)! + h(g(N) - 1)! > g(N)! + h(N - 1)! \ge g(N - 1)! + N - 1.
 \]
Therefore, by the definition of the sets $B_{i}$, we have $A_{i}(N) \le B_{i}(x) \le A_{i}(N+2h)$ for every $1 \le i \le h$. Thus we have,
 $B_{i}(x) = A_{i}(N) + O(1) = (1+o(1))A_{i}(N)$ as $x \rightarrow \infty$ for every $1 \le i \le h$. Now, we have $B_{1}(x) \cdots{} B_{h}(x) = (1 + o(1))A_{1}(N) \cdots{} A_{h}(N) = (1 + o(1))N$ as $x \rightarrow \infty$.
 It remains to prove that $B_{h+1}(x) = \frac{x}{N}(1 + o(1))$ as $x \rightarrow \infty$. It follows from the definition of $B_{h}$ that
 for $x \ge g(6)! + h(g(6)-1)!$ we have
 \[
 B_{h}(x) = g(6)! + h(g(6)-1)!
 \]
 \[
 + \sum_{n=6}^{N-1}\left(\frac{g(n+1)!+h(g(n+1)-1)!}{n-\lceil\sqrt{n}\rceil}-\frac{g(n)!}{n-\lceil\sqrt{n}\rceil}+3\right)
 \]
 \[
 + \left\lfloor\frac{x}{N-\lceil\sqrt{N}\rceil} - \frac{g(N)!}{N-\lceil\sqrt{N}\rceil}+3\right\rfloor
 \]
 \[
 = O(N) + \sum_{n=6}^{N-1}\left(\frac{g(n+1)!+h(g(n+1)-1)!}{n-\lceil\sqrt{n}\rceil}-\frac{g(n)!}{n-\lceil\sqrt{n}\rceil}\right)
 \]
 \[
 + \left(\frac{x}{N-\lceil\sqrt{N}\rceil} - \frac{g(N)!}{N-\lceil\sqrt{N}\rceil}\right).
 \]
 By $x \ge g(N)! + h(g(N)-1)! \ge N!$, we have $O(N) = o\left(\frac{x}{N}\right)$ as $x \rightarrow \infty$.
 It follows from (2) in Lemma 1 that $n \ge f_{h}(n)$. Then by the definition of $g(n)$, we have $g(n) \ge g(f_{h}(n) \ge n^{2}$. Applying this observation, a straightforward computation shows that
 \[
 \frac{g(n+1)!+h(g(n+1)-1)!}{n-\lceil\sqrt{n}\rceil}-\frac{g(n)!}{n-\lceil\sqrt{n}\rceil}
 \]
 \[
 = \left(1+O\left(\frac{1}{n^{2}}\right)\right)\cdot \frac{g(n+1)!+h(g(n+1)-1)!}{n-\lceil\sqrt{n}\rceil}
\]
\[
= \left(1+O\left(\frac{1}{\sqrt{n}}\right)\right)\cdot \frac{g(n+1)!}{n+1}.
\]
Hence,
\[
\sum_{n=6}^{N-1}\frac{g(n+1)!+h(g(n+1)-1)!}{n-\lceil\sqrt{n}\rceil}-\frac{g(n)!}{n-\lceil\sqrt{n}\rceil}
\]
\[
= \sum_{n = 6}^{N-1}\left(1+O\left(\frac{1}{\sqrt{n}}\right)\right)\cdot \frac{g(n+1)!}{n+1}.
\]
In the next step, we show that
\[
\sum_{n = 6}^{N-1}\left(1+O\left(\frac{1}{\sqrt{n}}\right)\right)\cdot \frac{g(n+1)!}{n+1} = (1 + o(1))\cdot \frac{g(N)!}{N}
\]
as $N \rightarrow \infty$.
Since $g(m)$ is strictly increasing,
\[
\frac{g(N+1)!}{g(N)!} \ge \frac{(g(N)+1)!}{g(N)!} = g(N) + 1 \ge N + 1 \ge \frac{N+1}{N},
\]
which implies that $\frac{g(N)!}{N}$ is monotone increasing. By $g(m) \ge m^{2}$, we have
\[
g(N-1)! \le \frac{1}{N^{2}}g(N)!.
\]
On the other hand,
\[
\frac{g(N-1)!}{N-1} \le \frac{g(N)!/N^{2}}{N-1} = O\left(\frac{g(N)!}{N^{3}}\right).
\]
By using the above observations, we have
\[
\sum_{n=6}^{N-1}\left(1+O\left(\frac{1}{\sqrt{n}}\right)\right)\cdot \frac{g(n+1)!}{n+1} = \sum_{n = 7}^{N-1}\left(1+O\left(\frac{1}{\sqrt{n}}\right)\right)\cdot \frac{g(n)!}{n} + \frac{g(N)!}{N}(1 + o(1))
\]
\[
= \sum_{n=7}^{N-1}O\left(\frac{g(N-1)!}{N-1}\right) + (1 + o(1))\frac{g(N)!}{N}
\]
\[
= O\left(N\frac{g(N)!}{N^{3}}\right) + \frac{g(N)!}{N}(1 + o(1)) = \frac{g(N)!}{N}(1 + o(1))
\]
as $x \rightarrow \infty$. It is clear that
\[
\frac{x}{N-\lceil\sqrt{N}\rceil} - \frac{g(N)!}{N-\lceil\sqrt{N}\rceil} = \left(1+O\left(\frac{1}{\sqrt{N}}\right)\right)\left(\frac{x-g(N)!}{N}\right)
\]
\[
= (1 + o(1))\frac{x-g(N)!}{N}
\]
as $x \rightarrow \infty$. Then it follows that
\[
B_{h+1}(x) = o\left(\frac{x}{N}\right) + (1 + o(1))\frac{g(N)!}{N} + \frac{x-g(N)!}{N}(1 + o(1)) = (1 + o(1))\frac{x}{N}
\]
as $x \rightarrow \infty$, which proves (iii). The proof of Lemma 1 is completed.

\subsection{Proof of Theorem 2}

Now, we prove Theorem 2 by induction on $h$. We show that there exist infinite sets
$A_{1}, \dots{} ,A_{h} \subset \mathbb{N}$ with the following properties:
\begin{itemize}
    \item[(1)] $A_{1} + \dots{} + A_{h} = \mathbb{N}$,
    \item[(2)] there exists a monotone increasing arithmetic function $f_{h}(n) \ge 0$ with $$\lim_{n\rightarrow \infty}f_{h}(n) = \infty$$ such that the equation
    $a_{1} + \dots{} + a_{h} = n$, $a_{i} \in A_{i}$ has a solution with
    $a_{i} \ge f_{h}(n)$.
    \item[(3)] $A_{1}(x) \cdots{}  A_{h}(x) = (1 + o(1))x$ as $x \rightarrow \infty$.
\end{itemize}
For $h = 1$ consider the set of natural numbers and the function $f_{1}(n) = n$, which gives the result.
Assume that the statement of Theorem 2 holds for $h$. For $h + 1$ the result follows from Lemma 1. (Actually, for $h = 2$ our construction is the same as the construction of Danzer \cite{Danzer}). The proof of Theorem 2 is completed.

\section{Proof of Theorem 4}
Let $h \ge 2$. We will prove that $\mathcal{A}_{h+1} \subseteq \mathcal{A}_{h}$.
Let $A \in \mathcal{A}_{h+1}$. Then there exist
$A_{2}, \dots{} ,A_{h+1} \subseteq \mathbb{N}$ such that $A + A_{2} + \dots{}  + A_{h+1} = \mathbb{N}$ and $A(x)A_{2}(x)\cdots{} A_{h+1}(x) = (1+o(1))x$ as $x \rightarrow \infty$. Let $A_{h}^{*} = A_{h} + A_{h+1}$. It is clear that $A^{*}_{h}(x) \le A_{h}(x) \cdot A_{h+1}(x)$. Then we have $A + A_{2} + \dots{}  + A_{h-1} + A_{h}^{*} = \mathbb{N}$ and so
$A(x)A_{2}(x)\cdots{}  A_{h-1}(x)A_{h}^{*}(x) \ge x + 1$. On the other hand,
$A(x)A_{2}(x)\cdots{}  A_{h-1}(x)A^{*}_{h}(x) \le A(x)A_{2}(x)\cdots{}  A_{h+1}(x) = (1 + o(1))x$ as $x \rightarrow \infty$, thus we have $A(x)A_{2}(x)\cdots{}  A_{h-1}(x)A_{h}^{*}(x) = (1 + o(1))x$ as $x \rightarrow \infty$, which implies that $A \in \mathcal{A}_{h}$. The proof of Theorem 4 is completed.


\begin{thebibliography}{99}
\bibitem{Chen} \textsc{Y.-G. Chen, J.-H. Fang}.\ \textit{On a conjecture of S\'ark\"ozy and Szemer\'edi}, Acta Arith., \textbf{169}, (2015) 47-58.
\bibitem{Danzer} \textsc{L. Danzer}.\ \textit{\"Uber eine Frage von G.Hanani aus der additiven Zahlentheorie}, J. Reine Angew. Math, \textbf{214/215}, (1964) 392-394.
\bibitem{erd} \textsc{P. Erd\H{o}s}.\ \textit{Some unsolved problems}, Michigan Math., J. \textbf{4} (1957) 291-300.
\bibitem{Nar} \textsc{W. Narkiewicz}.\ \textit{Remarks on a conjecture of Hanani in additive number theory}, Colloq. Math., \textbf{7}, (1959/60), 161-165.
\bibitem{Ruzsa2} \textsc{I. Z. Ruzsa}.\ \textit{An asymptotically exact additive completion}, Studia Sci. Math. Hungar., \textbf{32}, (1996) 51-57.
\bibitem{Ruzsa1} \textsc{I. Z. Ruzsa}.\ \textit{Additive completion of lacunary sequences}, Combinatorica, \textbf{21}, (2001) 279-291.
\bibitem{sarkoszeme} \textsc{A. S\'ark\"ozy, E. Szemer\'edi}.\ \textit{On a problem in additive number theory}, Acta Math. Hungar., \textbf{64}, (1994) 237-245.
\end{thebibliography}
\end{document}